%% file: 2013-SanLhoLem.tex
\documentclass[a4paper, final, natbib]{svjour3}

\usepackage[english]{babel}
\usepackage{hyperref}
\usepackage{graphicx}
\usepackage{enumerate}
\usepackage{amsmath, amsbsy, amsfonts}
\usepackage{graphics}
\usepackage{multicol}

\def\Zbb{\mathbb{Z}}
\def\Rbb{\mathbb{R}}
\def\Tbb{\mathbb{T}}
\def\Gscr{\mathcal{G}}
\def\Rscr{\mathcal{R}}
\def\phi{\varphi}
\def\rho{\varrho}

\begin{document}
\journalname{Celest Mech Dyn Astr (2014) 119:75–89}
\title{Effective stability around the Cassini state \\ 
in the spin-orbit problem}

\author{Marco Sansottera \and Christoph Lhotka \and Anne Lema{\^\i}tre}

\institute{M. Sansottera \and C.Lhotka \and A. Lema{\^\i}tre \at
  naXys, Department of Mathematics, University of Namur,\\
  Rempart de la Vierge, 8, 5000~Namur, Belgium \\
  Emails: marco.sansottera@unamur.be, lhotka@mat.uniroma2.it,
  anne.lemaitre@unamur.be
}

\maketitle

\begin{abstract}
We investigate the long-time stability in the neighborhood of the
Cassini state in the conservative spin-orbit problem.  Starting with
an expansion of the Hamiltonian in the canonical Andoyer-Delaunay
variables, we construct a high-order Birkhoff normal form and give an
estimate of the effective stability time in the Nekhoroshev sense.  By
extensively using algebraic manipulations on a computer, we explicitly
apply our method to the rotation of Titan.  We obtain physical bounds
of Titan's latitudinal and longitudinal librations, finding a stability
time greatly exceeding the estimated age of the Universe.  In
addition, we study the dependence of the effective stability time on
three relevant physical parameters: the orbital inclination, $i$, the
mean precession of the ascending node of Titan orbit, $\dot\Omega$,
and the polar moment of inertia, $C$.

\keywords{Spin-orbit resonance\and Normal form methods \and Cassini
  state \and Titan \and Long-time stability}

\end{abstract}

\section{Introduction}\label{sec:int}
Like the Moon, most of the regular satellites of the Solar System
present the same face to their
planet. Cassini~\citeyearpar{Cassini-1693}, considering a simplified
model of the rotation of the Moon, showed how this peculiar feature
corresponds to an equilibrium, called a \emph{Cassini state}.
Moreover, he showed how small perturbations of this model do not lead
to the destabilization of this equilibrium, but just to the excitation
of librations around it. The fact that so many natural satellites are
found in spin-orbit synchronization is due to dissipative effects,
see, e.g.,~Goldreich~\& Peale~\citeyearpar{GolPea-1966}.

Further investigations of the spin-orbit problem, also introducing
more refined models, have been carried out, see,
e.g.,~Colombo~\citeyearpar{Colombo-1966},
Peale~\citeyearpar{Peale-1969},
Beletskii~\citeyearpar{Beletskii-1972}, Ward~\citeyearpar{Ward-1975},
Henrard~\& Murigande~\citeyearpar{HenMur-1987}, Bouquillon
et~al.~\citeyearpar{BouKinSou-2003} and Henrard~\&
Schwanen~\citeyearpar{HenSch-2004}.  It has been found that the
equilibrium described by Cassini is not the only one.  Indeed,
depending on the values of some parameters, one or three other ones
are possible.

In this paper we study the problem of the stability around a Cassini
state in the light of Birkhoff normal form~\citeyearpar{Birkhoff-1927}
and Nekhoroshev theory~(\citeyear{Nekhoroshev-1977},
\citeyear{Nekhoroshev-1979}).  We remark that Nekhoroshev theorem
consists of an \emph{analytic part} and a \emph{geometric part\/}.  Our
approach is just based on the analytic part of the theorem;
nevertheless, let us stress that in the isochronous case one does not
need the geometric part to achieve the exponential stability time.

Our aim is to give a long-term analytical bound to both the
latitudinal and longitudinal librations.  Thus we give an estimate of
the \emph{effective stability time} for the elliptic equilibrium point
corresponding to the Cassini state.  The procedure that we propose is
reminiscent of the one that has already been used in Giorgilli
et~al.~(\citeyear{GioLocSan-2009}, \citeyear{GioLocSan-2010}) and
Sansottera et~al.~\citeyearpar{SanLocGio-2011} where the authors gave
an estimate of the long-time stability for the Sun-Jupiter-Saturn
system and the planar Sun-Jupiter-Saturn-Uranus system, respectively.
In addition, a study of the long-time stability for artificial
satellites has been done in Steichen~\&
Giorgilli~\citeyearpar{SteGio-1997}.

Our method has been applied to the largest moon of Saturn, Titan, by
making extensive use of algebraic manipulations on a computer.  As a
result, we obtain physical bounds for Titan latitudinal and
longitudinal librations.  It turns out that the effective stability
time around the Cassini state largely exceeds the estimated age of the
Universe.  Finally, we study the dependence of the stability time on
three relevant physical parameters: the orbital inclination, $i$, the
mean precession of the ascending node of Titan orbit, $\dot\Omega$,
and the polar moment of inertia, $C$.

The paper is organized as follows. In Section~\ref{sec:Hamiltonian},
we introduce the Hamiltonian model of the spin-orbit problem and
sketch the steps of the expansion in the canonical Andoyer-Delaunay
variables.  The Hamiltonian turns out to have the form of a perturbed
system of harmonic oscillators.  Therefore, we construct a high-order
Birkhoff normal form and investigate the stability around the Cassini
state in the Nekhoroshev sense.  This part is worked out in
Section~\ref{sec:stability}.  In Section~\ref{sec:Titan}, we apply our
method to Titan, giving an estimate of the effective stability time
and investigating its dependence on the mean inclination, precession
of the ascending node and polar moment of inertia.  Finally, our
results are summarized in Section~\ref{sec:results}.  An Appendix
follows.

\section{Hamiltonian formulation of the spin-orbit problem}\label{sec:Hamiltonian}
We consider a rotating body (e.g., Titan) with mass $m$ and equatorial
radius $R_{\rm e}\,$, orbiting around a point body (e.g., Saturn) with
mass $M$.  The rotating body is considered as a triaxial rigid body
whose principal moments of inertia are $A$, $B$ and $C$, with $A\leq
B<C$.

We closely follow the Hamiltonian formulation that has already been
used in previous works, see, e.g.,~Henrard~\&
Schwanen~\citeyearpar{HenSch-2004} for a general treatment of
synchronous satellites, Henrard~\citeyearpar{Henrard-2005a} for Io,
Henrard~\citeyearpar{Henrard-2005b} for Europa, D'Hoedt~\&
Lema{\^\i}tre~\citeyearpar{Dho04} and Lema{\^\i}tre
et~al.~\citeyearpar{LemDHoRam-2006} for Mercury, Noyelles
et~al.~\citeyearpar{NoyLemVie-2008} for Titan, and
Lhotka~\citeyearpar{Lhotka-2013} for a symplectic mapping model.  We
try to make the paper self-contained by describing the main aspects
and we refer to the quoted works for a more detailed exposition.
First, we briefly describe how to express the Hamiltonian in the
Andoyer-Delaunay set of coordinates.  Then we introduce a simplified
spin-orbit model that keeps the main features of the problem and
represents the basis of our study.

\subsection{Reference frames}\label{sbs:ReferenceFrames}
In order to describe the spin-orbit motion we need four reference
frames, centered at the center of mass of the rotating body, see
D'Hoedt~\& Lema{\^\i}tre~\citeyearpar{Dho04}:
\begin{enumerate}[(i)]
\item the \emph{inertial frame}, $(X_0, Y_0, Z_0)\,$, with $X_0$ and
  $Y_0$ in the ecliptic plane;
\item the \emph{orbital frame}, $(X_1, Y_1, Z_1)\,$, with $Z_1$
  perpendicular to the orbit plane;
\item the \emph{spin frame}, $(X_2, Y_2, Z_2)\,$, with $Z_2$ pointing to
  the spin axis direction and $X_2$ to the ascending node of the
  equatorial plane on the ecliptic plane;
\item the \emph{body frame}, $(X_3, Y_3, Z_3)\,$, with $Z_3$ in the
  direction of the axis of greatest inertia and $X_3$ of the axis of
  smallest inertia.
\end{enumerate}
We define the direction $\nu_{01}$ chosen along the ascending node of
the plane $(X_0,Y_0)$, on the plane $(X_1,Y_1)\,$.  Likewise we define
$\nu_{23}$ taking the planes $(X_2,Y_2)$ and $(X_3,Y_3)$.  We denote
by $K$ the angle between $Z_0$ and $Z_2$, and $J$ the angle between
$Z_2$ and $Z_3$.  We also introduce the angles related to the
rotational motion:~(i) $l_s$, between $\nu_{23}$ and $X_3$, measured
in the plane $(X_3, Y_3)$;~(ii) $g_s$, between $X_2$ and $\nu_{23}$,
measured in the plane $(X_2, Y_2)$;~(iii) $h_s$, between $X_0$ and
$X_2$, measured in the plane $(X_0, Y_0)$.  Finally, for the
orbital dynamics, we denote by $l_o$ the mean anomaly, $\omega$ the
pericenter argument and $\Omega$ the argument of the ascending node of
the rotating body.

The four reference frames and angles defined here above, that are
needed to introduce the Andoyer and Delaunay canonical variables, are
reported in Figure~\ref{fig:sop}.

\begin{figure}
\centering
\vskip5pt
\includegraphics[width=.35\linewidth]{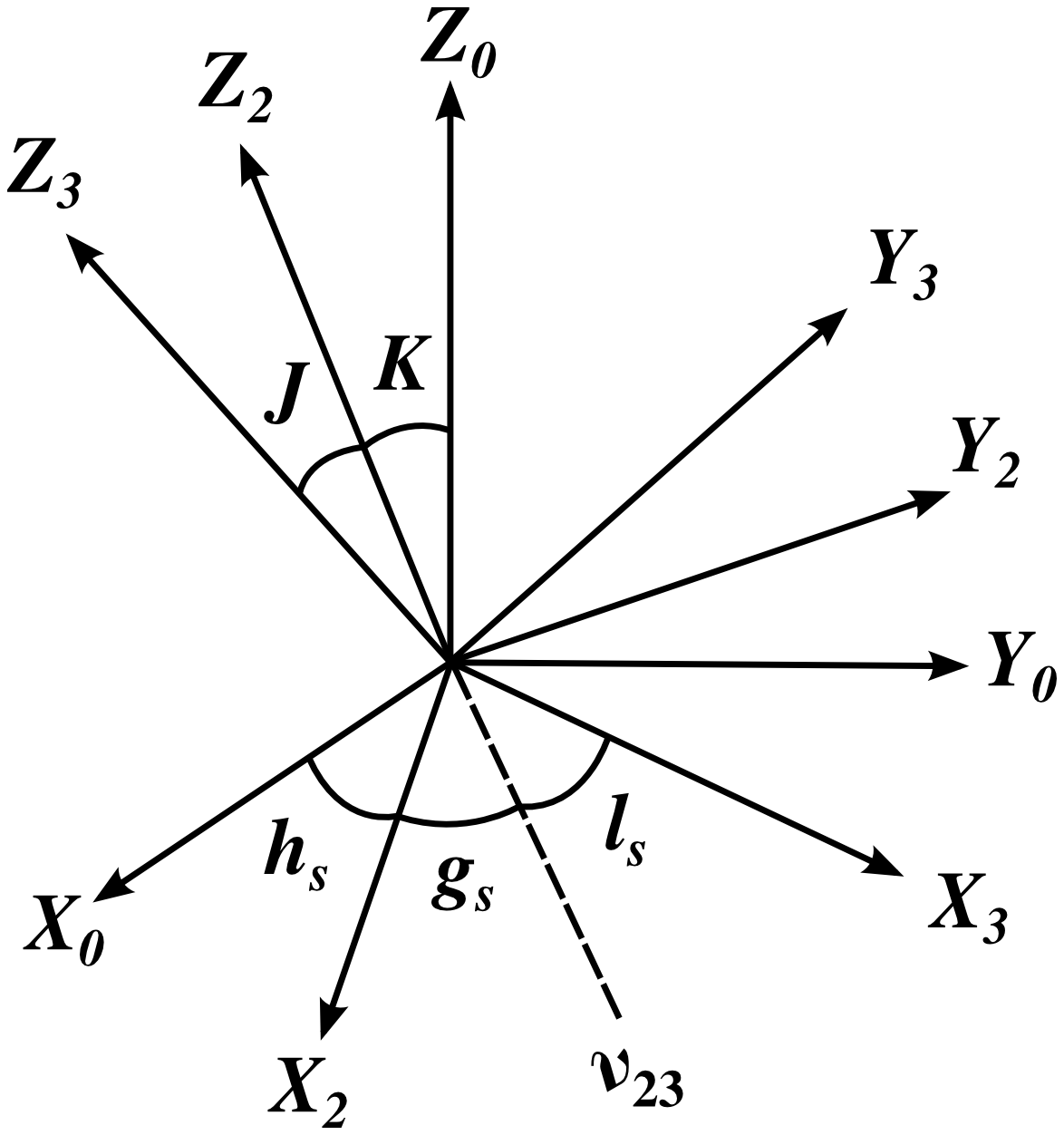}
\hspace{30pt}
\includegraphics[width=.35\linewidth]{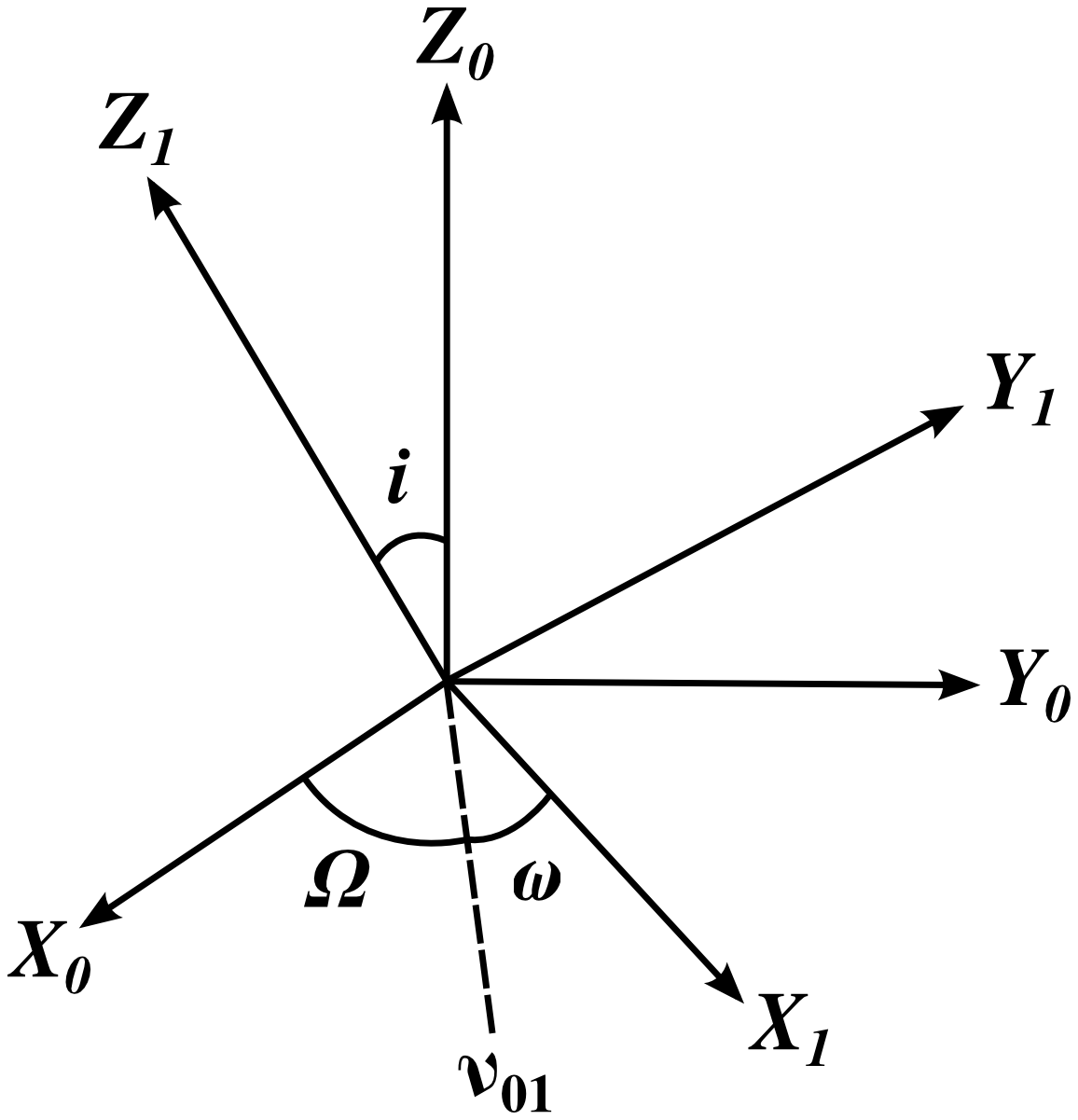}
\caption{The four reference frames and the relevant angles related to
  the Andoyer (left) and Delaunay (right) canonical variables.  See
  the text for more details.}\label{fig:sop}
\end{figure}

\subsection{Andoyer variables}
For the rotational motion, following Deprit~\citeyearpar{Deprit-1967}, we adopt the
Andoyer variables,
\begin{equation*}
\begin{aligned}
L_s &= G_s \cos{J}\ ,   &l_s\ ,\\
G_s &\ , &g_s\ ,\\
H_s &= G_s\cos{K}\ ,      &h_s\ ,
\end{aligned}
\end{equation*}
where $G_s$ is the norm of the angular momentum.

The Andoyer variables present two virtual singularities: for $J=0$,
i.e., when the angular momentum is collinear with the direction of the
axis of greatest inertia, and $K=0$, i.e., when the spin axis of
rotation is perpendicular to the orbital plane.  Therefore, we
introduce the modified Andoyer variables, see Henrard~\citeyearpar{Henrard-2005a},
\begin{equation}
\begin{aligned}
L_1 &= \frac{G_s}{n_o C}     \ , &l_1 &=  l_s+g_s+h_s\ ,\\
L_2 &= \frac{G_s-L_s}{n_o C} \ , &l_2 &= -l_s\ ,\\
L_3 &= \frac{G_s-H_s}{n_o C} \ , &l_3 &= -h_s\ ,\\
\end{aligned}
\label{eq:Andoyer}
\end{equation}
where $n_o$ is the orbital mean motion of the rotating body.  Let us
remark that this change of coordinates is a canonical transformation
with multiplier $(n_0 C)^{-1}$.

\subsection{Delaunay variables}
For the orbital motion, we introduce the classical Delaunay
variables,
\begin{equation*}
\begin{aligned}
L_o &= m\sqrt{\mu a}\ , &l_o&\ ,\\
G_o &= L_o\sqrt{1-e^2}\ , &g_o&=\omega\ ,\\
H_o &= G_o\cos{i}\ , &h_o&=\Omega\ ,
\end{aligned}
\end{equation*}
where $\mu=\Gscr(m+M)$, $\Gscr$ is the universal constant of
gravitation, $a$ is the semi-major axis, $e$ the eccentricity and $i$
the inclination of the orbit described by the rotating body.

The Delaunay variables are singular for $e=0$ and $i=0$.  Therefore,
we introduce the modified Delaunay variables,
\begin{equation}
\begin{aligned}
L_4 &= L_o     \ , &l_4 &=  l_o+g_o+h_o\ ,\\
L_5 &= L_o-G_o \ , &l_5 &= -g_o-h_o\ ,\\
L_6 &= G_o-H_o \ , &l_6 &= -h_o\ .
\end{aligned}
\label{eq:Delaunay}
\end{equation}
Let us remark that $L_5$ and $L_6$ are of the same order of magnitude
as the eccentricity, $e$, and inclination, $i$, respectively.

\subsection{Free rotation}
The rotational kinetic energy, Deprit~\citeyearpar{Deprit-1967}, reads
$$
T = \frac{L_s^2}{2C} +\frac{1}{2}
(G_s^2-L_s^2)
\left(
\frac{\sin^2{l_s}}{A}+\frac{\cos^2{l_s}}{B}
\right)\ .
$$
Thus, in our set of variables, we get
$$
\frac{T}{n_oC}= \frac{n_o L_1^2}{2} + \frac{n_0 L_3(2L_1-L_3)}{2}
\left( \frac{\gamma_1+\gamma_2}{1-\gamma_1-\gamma_2}\sin^2(l_3)
+ \frac{\gamma_1-\gamma_2}{1-\gamma_1+\gamma_2}\cos^2(l_3)\right)\ ,
$$
where
$$
\gamma_1=\frac{2C-B-A}{2C}
\quad\hbox{and}\quad
\gamma_2=\frac{B-A}{2C}\ .
$$

\subsection{Perturbation by another body}
The perturbation induced by the point body mass on the rotation of the
rigid body, can be expressed via a gravitational potential, $V$, see,
e.g.,~Henrard~\& Schwanen~\citeyearpar{HenSch-2004} for a detailed
exposition.  The gravitational potential has the form
$$
V = -{\Gscr M} \iiint_{W} \frac{\rho}{r}\,{\rm d}W\ ,
$$ where $\rho$ is the density inside the volume $W$ of the body and
$r$ the distance between the point mass and a volume element
inside the body.  Using the classical expansion of the potential in
spherical harmonics (see, e.g.,~Bertotti~\& Farinella\citeyearpar{BerFar-1990}), we get
$$
V = -\frac{\Gscr M}{r}\left(
1+\sum_{n\geq1}r^{-n}
\sum_{m=0}^{n} P_{n}^{m}(\sin\phi)\biggl(C_n^m \cos m\psi+S_n^m \sin m\psi\biggr)
\right)\ ,
$$
where $\phi$ is the latitude and $\psi$ the longitude of the
perturbing mass in the body frame.  Limiting the expansion to the
terms of order two and neglecting the first term $-{\Gscr M}/r$, which
does not affect the rotation, we have
$$
V = \frac{3}{2}\frac{\Gscr M}{a^3}\biggl(\frac{a}{r}\biggr)^3
\left(C_2^0 \bigl(x_3^2+y_3^2\bigr)-2 C_2^2 \bigl(x_3^2-y_3^2\bigr)
\right)\ ,
$$
where $(x_3,y_3,z_3)$ are the components (in the body frame) of
the unit vector pointing to the perturbing body.  The coefficients
$C_2^0$ and $C_2^2$, can be written in terms of the moments of inertia
and of the dimensionless parameters $J_2$ and $C_{22}$, as
\begin{align*}
C_2^0 &= \frac{A+B-2C}{2} = -m R_{{\rm e}}^2 J_2\ ,\\
C_2^2 &= \frac{B-A}{4} = m R_{{\rm e}}^2 C_{22}\ .
\end{align*}
Introducing the quantity $n_o^* = \sqrt{\Gscr M/a^3}$, we get
\begin{equation}
\frac{V}{n_o C} = n_0 \biggl(\frac{a}{r}\biggr)^3
\left(\delta_1 \bigl(x_3^2+y_3^2\bigr)+\delta_2 \bigl(x_3^2-y_3^2\bigr)
\right)\ ,
\label{eq:V_nC}
\end{equation}
with
$$
\delta_1=-\frac{3}{2} \left(\frac{n_o^*}{n_o} \right)^2\gamma_1
\quad\hbox{and}\quad
\delta_2=-\frac{3}{2} \left(\frac{n_o^*}{n_o} \right)^2\gamma_2\ .
$$
In order to express the potential $V$ in terms of the canonical
Andoyer-Delaunay variables, we use the following relation,
$$
\left(\begin{array}{c}x_3\\y_3\\z_3\end{array}\right)=
R_3(l_s) R_1(J) R_3(g_s) R_1(K) R_3(h_s) R_3(-h_o) R_1(-i_o) R_3(-g_o)
\left(\begin{array}{c}\cos(v_o)\\\sin(v_o)\\0\end{array}\right)\ ,
$$
where $v_o$ is the true anomaly and $R_i$ are rotation matrices,
$$
R_1(\phi) = 
\begin{pmatrix}
  \phantom{-}\cos(\phi) & \sin(\phi) & 0 \\
  -\sin(\phi) & \cos(\phi) & 0 \\
  0          & 0           & 1
 \end{pmatrix}
\qquad\hbox{and}\qquad
R_3(\phi) = 
\begin{pmatrix}
  1 &  0          & 0           \\
  0 &  \phantom{-}\cos(\phi) & \sin(\phi) \\
  0 & -\sin(\phi) & \cos(\phi)
 \end{pmatrix}\ .
$$ Finally, in order to get the expression of the potential in the
canonical variables, we have to compute the expansions of $(a/r)^3$,
$\cos(v_o)$ and $\sin(v_o)$ in terms of the eccentricity, $e$, and
mean anomaly, $l_o$, see, e.g. Poincar\'e~(\citeyear{Poincare-1892},
\citeyear{Poincare-1905}).

\subsection{The spin-orbit model}\label{sbs:model}
We now consider a simplified spin-orbit model, making some strong
assumptions on the system.  Similar assumptions have already been used
in previous studies, see, e.g.,~Henrard~\&
Schwanen~\citeyearpar{HenSch-2004}, Henrard~(\citeyear{Henrard-2005a},
\citeyear{Henrard-2005b}), Noyelles
et~al.~\citeyearpar{NoyLemVie-2008} and
Lhotka~\citeyearpar{Lhotka-2013}.

\begin{enumerate}[(i)]
\item We assume that the wobble, $J$, is equal to zero.  This means
  that the spin axis is aligned with respect to the direction of the
  axis of greatest inertia.  Thus, we have
$$
\frac{T}{n_o C} = \frac{n_o L_1^2}{2}\ .
$$
\item The rotating body is assumed to be locked in a $1$:$1$
  spin-orbit resonance.  We introduce the resonant variables
\begin{equation}
\begin{aligned}
\Sigma_1 &= L_1\ , &\sigma_1 &= l_1-l_4\ ,\\
\Sigma_3 &= L_3\ , &\sigma_3 &= l_3-l_6\ ,
\end{aligned}
\label{eq:resvar}
\end{equation}
and, after introducing these new coordinates, neglect the effects of
the fast dynamics on the long-term evolution via an average over
the angle, $l_4$, namely
$$
\langle V \rangle_{l_4} = \frac{1}{2\pi}\int_{0}^{2\pi}
V \,{\rm d}l_4\ .
$$
\item We neglect the influence of the rotation on the orbit of the
  body.  Indeed, we assume that the orbital variables, see
  Eq.~\eqref{eq:Delaunay}, are known functions of time acting as
  external parameters.  Namely, we consider that the perturbing body
  lies on a slowly precessing eccentric orbit, with constant
  precession frequency $\dot\Omega\,$.  The time dependence of the
  Hamiltonian can be modelled via the two angular variables,
  $$
  l_4(t) = n\,t + l_4(0)
  \qquad\hbox{and}\qquad
  l_6(t) = \dot\Omega\, t + l_6(0)\ .
  $$
\end{enumerate}

\noindent
Finally, we end up with a Hamiltonian that reads
\begin{equation}
H = \frac{n_o \Sigma_1^2}{2} -n_o \Sigma_1 + \dot\Omega \Sigma_3
 + \langle V \rangle_{l_4}\ .
\label{eq:H}
\end{equation}
Taking into account a more general perturbed orbit (e.g., including
the wobble and thus removing hypothesis~(i) by adding one degree of
freedom or replacing the first order average in~(ii) by an higher
order average), would not be difficult, but would introduce many
supplementary parameters.  The expanded Hamiltonian would have a form
similar to Eq.~\eqref{eq:H}, but, having more terms, it would be
heavier from the computational point of view.  In this work we want to
focus on the estimate of the long-time stability around the Cassini
state, thus we take this simplified spin-orbit model that keeps the
relevant features of the system.

The Hamiltonian~\eqref{eq:H} possesses an equilibrium, the Cassini
state, defined by
\begin{equation}
\begin{aligned}
\sigma_1 &= 0\ , &\frac{\partial H}{\partial \Sigma_1} &= 0\ ,\\
\sigma_3 &= 0\ , &\frac{\partial H}{\partial \Sigma_3} &= 0\ .
\end{aligned}
\label{eq:cassini}
\end{equation}
We denote by $\Sigma_1^*$ and $\Sigma_3^*$ the values at the
equilibrium.  Geometrically, $\sigma_1=0$ means that the smallest axis
of inertia points towards the perturbing body, while $\sigma_3=0$
means that the nodes of the orbit and equator are locked.  The values
$\Sigma_1^*$ and $\Sigma_3^*$, correspond to fixing the inertial
obliquity $K^*$.  Finally, $\Sigma_1^*-1$ is a small correction of the
unperturbed spin frequency, $n_0 \Sigma_1^*$.

\section{Stability around the Cassini state}\label{sec:stability}
We now aim to study the dynamics in the neighborhood of the Cassini
state defined here above.  We introduce the translated canonical
variables
\begin{equation}
\begin{aligned}
&\Delta\Sigma_1 = \Sigma_1-\Sigma_1^* \ , \qquad \sigma_1\ ,\\
&\Delta\Sigma_3 = \Sigma_3-\Sigma_3^* \ , \qquad \sigma_3\ ,
\end{aligned}
\label{eq:cassini_shift}
\end{equation}
and, with a little abuse of notation, in the following we will denote
again $\Delta\Sigma_i$ by $\Sigma_i$, with $i=1,\,3\,$.  Let us also
introduce the shorthand notations $\Sigma=(\Sigma_1, \Sigma_3)$ and
$\sigma=(\sigma_1, \sigma_3)\,$.  In these new coordinates, the
equilibrium is set at the origin, thus we can expand the
Hamiltonian~\eqref{eq:H} in power series of $(\Sigma,\,\sigma)$.  Let
us remark that the linear terms disappear, as the origin is an
equilibrium, thus the lowest order terms in the expansion are quadratic
in $(\Sigma,\,\sigma)\,$. Precisely, we can write the Hamiltonian as
\begin{equation}
H(\Sigma,\sigma) = H_0(\Sigma,\sigma) + \sum_{j>0} H_j(\Sigma,\sigma)\ ,
\label{eq:Hsigma}
\end{equation}
where $H_j$ is an homogeneous polynomial of degree $j+2$ in
$(\Sigma,\,\sigma)\,$.  In the latter equation the quadratic term,
$H_0$, has been separated in view of its relevance in the perturbative
scheme.  The analytical expression of $H_0$ can be found in the
Appendix A of Henrard~\& Schwanen~\citeyearpar{HenSch-2004}.

\subsection{Diagonalization of the quadratic part}\label{sbs:diagonalization}
The quadratic part of the Hamiltonian reads
\begin{align*}
H_0(\Sigma,\sigma) =
&\mu_{{\scriptscriptstyle \Sigma_1 \Sigma_1}} \Sigma_1^2
+2\mu_{{\scriptscriptstyle \Sigma_1 \Sigma_3}} \Sigma_1 \Sigma_3
+\mu_{{\scriptscriptstyle \Sigma_3 \Sigma_3}} \Sigma_3^2\\
&+\mu_{{\scriptscriptstyle \sigma_1 \sigma_1}} \sigma_1^2
+2\mu_{{\scriptscriptstyle \sigma_1 \sigma_3}} \sigma_1 \sigma_3
+\mu_{{\scriptscriptstyle \sigma_3 \sigma_3}} \sigma_3^2\ .
\end{align*}
Following the approach of Henrard~\&
Schwanen~\citeyearpar{HenSch-2004}, we introduce a canonical
transformation to reduce the Hamiltonian $H_0$ to a sum of squares,
\begin{equation}
\begin{aligned}
\Sigma_1 &= (1-\alpha\beta)\Sigma_1' - \alpha\Sigma_3'\ , &\sigma_1 &= \sigma_1' - \beta\sigma_3'\ ,\\
\Sigma_3 &= \beta\Sigma_1' +\Sigma_3'\ , &\sigma_3 &= \alpha\sigma_1' +(1-\alpha\beta)\sigma_3'\ ,
\end{aligned}
\label{eq:diagvar}
\end{equation}
where the parameter $\alpha$ and $\beta$ have to be
chosen\footnote{The analytical form of $\alpha$ and $\beta$ can be
  found in the Henrard~\& Schwanen~\citeyearpar{HenSch-2004} (see
  equations~(16) and~(17)).} to get rid of the mixed terms in $H_0$,
namely,
\begin{align*}
\beta\mu_{{\scriptscriptstyle \sigma_1 \sigma_1}}
-(1-2\alpha\beta)\mu_{{\scriptscriptstyle \sigma_1 \sigma_3}}
-\alpha(1-\alpha\beta)\mu_{{\scriptscriptstyle \sigma_3 \sigma_3}}=0\ ,\\
\alpha(1-\alpha\beta)\mu_{{\scriptscriptstyle \Sigma_1 \Sigma_1}}
-(1-2\alpha\beta)\mu_{{\scriptscriptstyle \Sigma_1 \Sigma_3}}
-\beta\mu_{{\scriptscriptstyle \Sigma_3 \sigma_3}}=0\ .
\end{align*}
This change of coordinates is the so-called ``untangling
transformation'', see Henrard~\&
Lema{\^\i}tre~\citeyearpar{HenLem2005}, and permits to write $H_0$ as
\begin{equation}
H_0(\Sigma',\sigma') =
\left(\mu'_{{\scriptscriptstyle \Sigma_1' \Sigma_1'}} \Sigma_1'^2
+\mu'_{{\scriptscriptstyle \sigma_1' \sigma_1'}} \sigma_1'^2\right)
+\left(\mu'_{{\scriptscriptstyle \Sigma_3' \Sigma_3'}} \Sigma_3'^2
+\mu'_{{\scriptscriptstyle \sigma_3' \sigma_3'}} \sigma_3'^2\right)\ .
\label{eq:untangling}
\end{equation}
Thus, if the products of the coefficients $\mu'_{{\scriptscriptstyle
    \Sigma_i' \Sigma_i'}}$ and $\mu'_{{\scriptscriptstyle \sigma_i'
    \sigma_i'}}$, with $i=1,\,3\,$, are positive the quadratic part
of the Hamiltonian describes a couple of harmonic oscillators.

We now perform a rescaling and introduce the polar coordinates,
\begin{equation}
\begin{aligned}
\Sigma_1' &= \sqrt{2U_1/U_1^*} \cos(u_1)\ , \qquad\sigma_1' = \sqrt{2 U_1 U_1^*} \sin(u_1)\ ,\\
\Sigma_3' &= \sqrt{2U_3/U_3^*} \cos(u_3)\ , \qquad\sigma_3' = \sqrt{2 U_3 U_3^*} \sin(u_3)\ ,
\end{aligned}
\label{eq:actang}
\end{equation}
where
$$
U_1^* = \sqrt{\mu'_{{\scriptscriptstyle \Sigma_1' \Sigma_1'}}/\mu'_{{\scriptscriptstyle \sigma_1' \sigma_1'}}}\ ,
\qquad\hbox{and}\qquad
U_3^* = \sqrt{\mu'_{{\scriptscriptstyle \Sigma_3' \Sigma_3'}}/\mu'_{{\scriptscriptstyle \sigma_3' \sigma_3'}}}\ .
$$
After this last transformation, the quadratic part of the
Hamiltonian is expressed in action-angle variables and reads
$$
H_0 = \omega_{u_1} U_1 + \omega_{u_3} U_3\ ,
$$ where $\omega_{u_1}$ and $\omega_{u_3}$ are the frequencies of the
angular variables $u_1$ and $u_3$, respectively.  Again, we introduce
the shorthand notations $U=(U_1,\,U_3)\,$, $u=(u_1,\,u_3)$ and
$\omega_{u}=(\omega_{u_1},\,\omega_{u_3})\,$.

Finally, we apply the same transformations that allow to introduce
the action-angle variables for the quadratic term, i.e.,
equations~\eqref{eq:untangling} and~\eqref{eq:actang}, to the full
Hamiltonian~\eqref{eq:Hsigma}.  In these new coordinates,
the transformed Hamiltonian can be expanded in Taylor-Fourier series
and reads
\begin{equation}
H^{(0)}(U,\,u) = \omega_{u}\cdot U
+ \sum_{j>0} H^{(0)}_{j}(U,\,u)\ ,
\label{eq:H_Uu}
\end{equation}
where the terms $H_j$ are homogeneous polynomials of degree $j/2+1$ in
$U$, whose coefficients are trigonometric polynomials in the angles
$u\,$.

The Hamiltonian~\eqref{eq:H_Uu} has the form of a perturbed system of
harmonic oscillators, thus we are led to study the stability of the
equilibrium, placed at the origin, corresponding to the Cassini state.
Let us stress that the change of coordinates in Eq.~\eqref{eq:actang} is
singular at the origin, nevertheless, this virtual singularity is
harmless as we will just be interested in giving a bound of an
analytic function in a disc around the origin.  Moreover, one could
adopt the cartesian coordinates avoiding the singularity problem.  We
now aim to investigate the stability of the equilibrium in the light
of Nekhoroshev theory, introducing the so-called \emph{effective
  stability time}.

\subsection{Birkhoff normal form}
Following a quite standard procedure we construct the Birkhoff normal
form for the Hamiltonian~\eqref{eq:H_Uu}
(see~Birkhoff~\citeyearpar{Birkhoff-1927}; for an application of
Nekhoroshev theory see, e.g.,~Giorgilli~\citeyearpar{Giorgilli-1988}).
This is a well known topic, thus we limit our exposition to a short
sketch adapted to the present context.

The aim is to give the Hamiltonian the normal form at order $r$
\begin{equation}
H^{(r)}(U,\,u) = Z_0(U)
+\ldots+Z_r(U)+\sum_{s>r} \Rscr^{(r)}_{s}(U,\,u)\ ,
\label{eq:H_r}
\end{equation}
where $Z_s$, for $s=0,\,\ldots\,,r\,$, is a homogeneous polynomial of
degree $s/2+1$ in $U$ and in particular it is zero for odd $s$. The
unnormalized reminder terms $\Rscr^{(r)}_s$, for $s>r$, are
homogeneous polynomials of degree $s/2+1$ in $U$, whose coefficients
are trigonometric polynomials in the angles $u\,$.

We proceed by induction. For $r=0$ the Hamiltonian~\eqref{eq:H_Uu} is
already in Birkhoff normal form.  Assume that the Hamiltonian is in
normal form up to a given order $r>0$, we determine the generating
function $\chi^{(r+1)}$ and the normal form term $Z_{r+1}$, by solving
the equation
$$
\left\{ \chi^{(r+1)},\, \omega_{u}\cdot U \right\} +
 \Rscr^{(r)}_{r+1}(U,\,u) = Z_{r+1}(U)\ ,
$$
where $\{\cdot,\cdot\}$ denotes the usual Poisson bracket.  Using the
Lie series algorithm, see, e.g.,~Henrard~\citeyearpar{Henrard-1973} and
Giorgilli~\citeyearpar{Giorgilli-1995}, we compute the new Hamiltonian as $H^{(r+1)} =
\exp(L_{\chi^{(r+1)}}) H^{(r)}\,$.  It is easy to show that
$H^{(r+1)}$ has a form analogous to Eq.~\eqref{eq:H_r} with new functions
$\Rscr^{(r+1)}_s$ of degree $s/2+1$ (with $s > r+1$) and the normal
form part ending with $Z_{r+1}$ , which is equal to zero if $r$ is
even.  As usual when using Lie series methods, with a little abuse
of notation, we denote again by $(U,\,u)$ the new coordinates.

This algorithm can be iterated up to the order $r$ provided that the
non-resonance condition
$$
k\cdot\omega_u\neq0 \qquad\hbox{for } k=(k_1,\,k_2)\in \Zbb^2\hbox{ such that } 0<|k|_1\leq r+2\ ,
$$
is fulfilled, where we used the usual notation $|k|_1=|k_1|+|k_2|$.

\subsection{Effective stability}\label{sbs:time}
It is well known that the Birkhoff normal form at any finite order $r$
is convergent in some neighborhood of the origin, but the analyticity
radius shrinks to zero when the order $r\to\infty\,$.  Thus, the best
strategy is to look for stability over a finite time, possibly long
enough with respect to the lifetime of the system.  We concentrate
here on the quantitative estimates that allow to give an upper bound
of the effective stability time that has to be evaluated.

Let us pick two positive numbers $R_1$ and $R_3$, and consider a
polydisk $\Delta_{\rho R}$ centered at the origin of $\Rbb^2$, defined
as
$$
\Delta_{\rho R} = \left\{
U\in \Rbb^2 : |U_j| \leq \rho R_j\ ,\ j=1,3
\right\}\ ,
$$
$\rho>0$ being a parameter.

Let us consider a function
$$
f_s(U,\,u)=\sum_{|l|=s+2,k\in\Zbb^2} f_{l,k}
\,U^{l/2} {\sin\atop\cos} (k\cdot u)\ ,
$$ which is a homogeneous polynomial of degree $s/2+1$ in
the actions $U$ and depends on the angles $u$.  We define the quantity
$|f_s|_R$ as
$$
|f_s|_R = \sum_{|l|=s+2,k\in\Zbb^2} |f_{l,k}| R_1^{l_1/2} R_3^{l_2/2}\ .
$$
Thus we get the estimate
$$
|f(U,u)| \leq |f|_R\, \rho^{s/2+1}\ ,
\quad\hbox{for }
U\in\Delta_{\rho R}\,,\ u\in\Tbb^2\ .
$$

Let now $U(0)\in\Delta_{\rho_0 R}$, with $\rho_0 < \rho$.  Then we
have $U(t)\in\Delta_{\rho R}$ for $t\leq T$, where $T$ is the escape
time from the domain $\Delta_{\rho R}$.  This is the effective
stability time that we want to evaluate.  We consider the trivial
estimate
$$
|U(t)-U(0)|\leq |t|\cdot\sup|\dot U|\ ,
$$
thus we need to bound the quantity $\sup|\dot U|$.  To this end,
taking the Hamiltonian~\eqref{eq:H_r}, which is in Birkhoff
normal form up to order $r$, we get
$$
|\dot U| \leq |\{ U, H^{(r)} \}| = \sum_{s>r} |\{ U, \Rscr^{(r)}_s \}|
\leq c |\{ U, \Rscr^{(r)}_{r+1} \}|_R\, \rho^{r/2+1}\ ,
$$
with $c\geq1$.  In fact, after having set $\rho$ smaller than the
convergence radius of the remainder series, $\Rscr^{(r)}_s$ for $s>r$,
the above inequality holds true for some value $c$.

The latter equation allows us to find a lower bound for the escape
time from the domain $\Delta_{\rho R}$, namely
\begin{equation}
\tau(\rho_0, \rho, r) = 
\frac{\rho-\rho_0}{c |\{ U, \Rscr^{(r)}_{r+1} \}|_R\, \rho^{r/2+1}}\ ,
\label{eq:time}
\end{equation}
which, however, depends on $\rho_0$, $\rho$, and $r$.  Let us
emphasize that $\rho_0$ is the only physical parameter, being fixed by
the initial data, while $\rho$, and $r$ are left arbitrary.  Indeed,
the parameter $\rho_0$ must be chosen in such a way that the domain
$\Delta_{\rho_0 R}$ contains the initial conditions of the system.
Thus we try to find an estimate of the escape time, $T(\rho_0)$,
depending only on the physical parameter $\rho_0$.

We optimize $\tau(\rho_0, \rho, r)$ with respect to $\rho$ and $r$,
proceeding as follows. First we keep $r$ fixed, and remark that the
function $\tau(\rho_0, \rho, r)$ has a maximum for
$$
\rho_{\rm opt} = \frac{r+2}{r} \rho_0\ .
$$
This gives an optimal value of $\rho$ as a function of $\rho_0$ and
$r$, thus we define
$$
\widetilde\tau(\rho_0, r) = \tau(\rho_0, \rho_{\rm opt}, r)\ .
$$

Finally, we look for the optimal value $r_{\rm opt}$ of the
normalization order, which maximizes $\widetilde\tau(\rho_0, r)$.
Namely, we look for the quantity
$$
T(\rho_0) = \max_{r\geq1} \widetilde\tau(\rho_0, r)\ ,
$$
which is our best estimate of the escape time, we define this quantity
as the \emph{effective stability time}.

\section{Application to Titan}\label{sec:Titan}
We now come to the application of our study to the largest moon of
Saturn, Titan.  Let us stress that we take a simplified model: we
consider Titan as a rigid body orbiting around Saturn, regarded as a
point body mass.  Moreover, we keep all the assumptions that we have made
in Subsection~\ref{sbs:model} for the spin-orbit model.  Concerning
the hypothesis~(iii), i.e., the constant precession frequency
$\dot\Omega\,$, we remark that in this case is
not restrictive.  Indeed, for Titan, the expansion of
$\dot\Omega$ has a dominant term, while the other ones are negligible,
see~Vienne~\& Duriez~\citeyearpar{VieDur-1995} (Table~6d).

We first make an expansion of the Hamiltonian~\eqref{eq:H} up to
degree $8$ in the eccentricity, by using the Wolfram Mathematica
software.  Then, we express the Hamiltonian in the canonical
Andoyer-Delaunay variables, see equations~\eqref{eq:Andoyer}
and~\eqref{eq:Delaunay}, and introduce the resonant coordinates
$(\Sigma,\,\sigma)$, defined in Eq.~\eqref{eq:resvar}.  In the actual
computation we take the physical parameters reported in
Table~\ref{tab:Titan_param}.

\renewcommand{\arraystretch}{1.2}
\begin{table}[!h]
\begin{center}
\caption{Titan physical parameters.  The parameters related to the
  gravity field and shape of Titan, $J_2$, $C_{22}$ and $C/MR_{\rm
    e}^2$, are taken from Iess
  et~al.~\citeyearpar{IesRapJacRacSte-2010}.  The orbital data are
  taken from Vienne~\& Duriez~\citeyearpar{VieDur-1995}.}
\begin{tabular}{|c|l|}
\hline
\vphantom{$|^{|^|}$} $M$ & $\phantom{-}5.683200\times10^{+26}\,$ Kg \\
$J_2$ & $\phantom{-}3.180800\times10^{-5}$ \\
$C_{22}$ & $\phantom{-}9.983000\times10^{-6}$ \\
$C/MR_{\rm e}^2$ & $\phantom{-}3.414000\times10^{-1}$ \\
$a$ & $\phantom{-}1.221865\times10^{+6}\,$ Km\\
$e$ & $\phantom{-}2.890000\times10^{-2}$ \\
$i$ & $\phantom{-}5.579818\times10^{-3}\,$ rad\\
$\dot\Omega$ & $-8.931240\times10^{-3}\,$ rad/year\\
$n_o$ & $\phantom{-}1.439240\times10^{+2}\,$ rad/year\\
\hline
\end{tabular}
\label{tab:Titan_param}
\end{center}
\end{table}

A numerical solution of Eq.~\eqref{eq:cassini} gives $\Sigma_1^* =
1+1.829\times10^{-9}$ and $\Sigma_3^* = 2.947\times10^{-5}$.  These
are the values at the equilibrium, the Cassini state; thus we expand
the Hamiltonian~\eqref{eq:Hsigma} in the translated canonical
variables, see Eq.~\eqref{eq:cassini_shift}.  Following the procedure in
Subsection~\ref{sbs:diagonalization}, we perform the so-called
``untangling transformation'', see Eq.~\eqref{eq:diagvar} and
introduce the action-angle variables $(U,\,u)\,$, where
$U_1^*=5.348\times10^1$ and $U_3^*=1.696\times10^4$.  This allows us
to write the Hamiltonian as
\begin{equation*}
H^{(0)}(U,\,u) = \omega_{u}\cdot U
+ \sum_{j>0} H^{(0)}_{j}(U,\,u)\ ,
\end{equation*}
with $\omega_{u_1} = 2.690$ and $\omega_{u_3} = 2.375\times10^{-2}\,$.
More details concerning the semi-analytical expansion of the
Hamiltonian, including the explicit form of some relevant quantities,
are reported in Appendix~\ref{sec:A.1}.  We now compute a
high-order Birkhoff normal form up to order $r=30$,
see Eq.~\eqref{eq:H_r}, by using a specially devised algebraic
manipulator, see Giorgilli~\& Sansottera~\citeyearpar{GioSan-2012}.
As shown in Subsection~\ref{sbs:time}, the estimate of the effective
stability time is now straightforward.  Let us remark that,
in this specific case, it is enough to set $c=2$ in Eq.~\eqref{eq:time}.

First, we give an estimate of the stability time as a function of the
parameter $\rho_0$, that parametrizes the radius of the polydisk
containing the initial data.  Let us remark that $\rho_0 = 0$
corresponds to the exact Cassini state, while taking $\rho_0 > 0$
allows small oscillations around the equilibrium point.  The results
of our computations are reported in Figure~\ref{fig:time}.  The best
estimate corresponds to the highest order of normalization, namely
$r=30$, however already at order $r=20$ we obtain really good
estimates.  As shown in the plot, we reach an effective stability time
greatly exceeding the estimated age of the Universe in a domain
$\Delta_{R}$ that roughly corresponds to a libration of $0.1$ radians.
Quoting J.E. Littlewood in his papers about the stability of the
Lagrangian equilateral equilibria of the problem of three bodies,
\emph{``while not eternity, this is a considerable slice of it.''}.

\begin{figure}
\begin{center}
\resizebox{.8\textwidth}{!}{\input{time_rho0_U3.tex}}
\renewcommand\figurename{}
\renewcommand\thefigure{2(a)}
\caption{}
\label{fig:time}

\resizebox{.9\textwidth}{!}{\input{time_Omdot-incl_U3.tex}}
\vskip-10pt 
\renewcommand\figurename{}
\renewcommand\thefigure{2(b)}
\caption{}
\label{fig:timedep}
\end{center}
\setcounter{figure}{1}
\caption{Evaluation of the estimated effective stability time.  The
  time unit is the year.  In~\ref{fig:time} we report the effective
  stability time evaluated according to the algorithm described in
  Subsection~\ref{sbs:time}. The three lines correspond (from down to
  top) to three different normalization orders: $r=10$ (blue), $r=20$
  (pink) and $r=30$ (red).  In~\ref{fig:timedep}--\ref{fig:timedep3}
  we report the effective stability time as a function of the mean
  inclination, $i$, (\textit{continued on the following page})}
\end{figure}

\begin{figure}
\begin{center}
\resizebox{.9\textwidth}{!}{\input{time_Omega-CRM_U3.tex}}
\vskip-10pt
\renewcommand\figurename{}
\renewcommand\thefigure{2(c)}
\caption{}
\label{fig:timedep2}
\resizebox{.9\textwidth}{!}{\input{time_CRM-incl_U3.tex}}
\vskip-10pt
\renewcommand\figurename{}
\renewcommand\thefigure{2(d)}
\caption{}
\label{fig:timedep3}
\end{center}
\setcounter{figure}{1}
\caption{(\textit{continued from the previous page}) the mean
  precession of the ascending node of Titan orbit, $\dot\Omega\,$ and
  the normalized greatest moment of inertia, $C/MR_{\rm e}^2$.  The
  colors refer to the stability time (in logarithmic scale), the
  yellow region is the most stable one.}
\end{figure}

\setcounter{figure}{3}

Finally, we investigate the dependence of the effective stability time
on three relevant physical parameters: the orbital inclination, $i$,
the mean precession of the ascending node of Titan orbit,
$\dot\Omega$, and the polar moment of inertia, $C$.

The results for the investigation of the stability time as a function
of $i$ and $\dot\Omega\,$, is reported in Figure~\ref{fig:timedep}.
The effective stability time turns out to grow up for increasing
values of the inclination, while the dependence on the precession of
the ascending node is less significant.  In Figure~\ref{fig:timedep2}
we report the outcome of our computations of the effective stability
time as a function of $C/MR_{\rm e}^2$ and $\dot\Omega\,$.  We let
vary the parameter $C/MR_{\rm e}^2$ in the range $[0.3,0.4]$.  This
interval should guarantee that we include the \emph{true} value of the
parameter, as for the Galilean satellites of Jupiter.  We see that the
effective stability time decreases for both increasing values of the
normalized greatest moment of inertia and precession of the ascending
node.  The dependence of the stability time on the parameters $i$
and $C/MR_{\rm e}^2$ is reported in Figure~\ref{fig:timedep3}.  Again
the effective stability time is growing for increasing values of the
inclination, while the parameter $C/MR_{\rm e}^2$, in this case, is
less important.

\section{Conclusions and outlook}\label{sec:results}
We have studied the long-time stability around the Cassini state in
the spin-orbit problem.  As explained in details in
Section~\ref{sec:Hamiltonian} and~\ref{sec:stability}, the Hamiltonian
of the system turns out to have the form of a perturbed system of
harmonic oscillators.  We have computed a high-order Birkhoff normal
form that allows us to obtain an analytical estimate of the effective
stability time around the Cassini state.

Our aim was to investigate the physical relevance of the long-time
stability in the framework of the Nekhoroshev theory.  For this reason
we used a simplified spin-orbit model that keeps all the relevant
features of the system.  Taking a more general model into account,
would not substantially change the form of the Hamiltonian and
certainly deserves to be studied in the future.

The main conclusion is that our estimates are physically relevant.  We
investigate a region of the parameters that contains the real data of
Titan and the stability around the Cassini state is assured over very
long times, largely exceeding the age of the Universe.  Therefore, our
work is a further confirmation that the long-time stability predicted
by Nekhoroshev theory may be very relevant for physical systems.

The effective stability time depends on many physical parameters,
e.g., the mean inclination, $i$, the mean precession rate,
$\dot\Omega$ and the greatest moment of inertia, $C$.  These
parameters play a crucial role in understanding Titan rotation and
their estimation is a major challenge.  As shown in the last part of
the previous section, our method can be used to determine the most
stable region of the parameters and supports the estimates given by
observations.

The natural question is how far these results can be extended
considering more realistic models.  On the one hand, let us recall
that the simplified spin-orbit model we adopt, makes strong
assumptions on the system.  Thus, our model might be less perturbed
than the \emph{real} system, and therefore more stable.  However, if
the question concerns the stability around the Cassini state, meaning
that we just look for bounds on the librations, then further terms in
the Hamiltonian would be relevant only if they produce resonances that
make their size in the perturbation expansions to grow, due to the
small divisors.  Indeed, all our estimates of the stability time are
based on the evaluation of the size of the perturbing terms.  Finally,
let us remark that our method is strongly based on Hamiltonian theory,
thus other non conservative effects cannot be taken into account in
this framework.

\begin{acknowledgements}
The work of C.~L. was financially supported by the contract Prodex
C90253 ``ROMEO'' from BELSPO, and partly by the Austrian FWF research
grant P-J3206.  The work of M.~S. is supported by an FSR Incoming
Post-doctoral Fellowship of the Acad\'emie universitaire Louvain,
co-funded by the Marie Curie Actions of the European Commission.
\end{acknowledgements}

\appendix

\section{Semi-analytical expansion of the Hamiltonian (Titan application)}\label{sec:A.1}
We report here the explicit expansions of the relevant Hamiltonian
functions related to the application to Titan.  We recall that the
Titan physical parameters adopted here are reported in
Table~\ref{tab:Titan_param}.

The averaged potential in Eq.~\eqref{eq:H} takes the form
\begin{align*}
\langle V \rangle_{l_4}=&-1.01\times 10^{-2}\cos^2 K-3.13\times 10^{-7} \sin^2 K \\
 &-1.12\times 10^{-4} \cos \left(\sigma _3\right) \cos K \sin K \\
 &-\cos \left(2 \sigma_1\right)\left(3.14\times10^{-3}+6.29\times10^{-3}\cos K+3.14\times10^{-3}\cos^2 K\right)\\
 &+\cos\left(2 \sigma _3\right) \left(1.56\times 10^{-7}\cos^2 K-1.56\times10^{-7}\right) \\
 &-\cos \left(2 \sigma _1+\sigma _3\right)\left(3.51\times 10^{-5}\cos K+3.51\times 10^{-5}\right) \sin K \\
 &+\cos\left(2 \sigma _1+2 \sigma _3\right)\left(4.90\times10^{-8}\cos^2 K-9.79\times10^{-8}\sin^2 K-4.90\times10^{-8}\right) \\
 &+\cos \left(2 \sigma _1+3\sigma _3\right) \left(2.73\times10^{-10} \cos K-2.73\times 10^{-10}\right)\sin K \\
 &-1.91\times 10^{-13} \cos \left(2 \sigma _1+4 \sigma_3\right)\left(1.00-1.00 \cos K\right)^2\ .
\end{align*}

The quadratic part of the Hamiltonian~\eqref{eq:Hsigma}, $H_0$, reads
\begin{align*}
H_0 &= 2.52\times 10^{-2}\sigma_1^2 + 1.08\times 10^{-6} \sigma_1\sigma_3 + 7.00\times 10^{-7} \sigma_3^2 \\
 &\quad + 7.20\times 10^1 \Sigma_1^2  - 2.08\times 10^{-2}\Sigma_1 \Sigma_3 + 2.01\times 10^2\Sigma_3^2\ ,
\end{align*}
and the values of the parameter $\alpha$ and $\beta$ corresponding to
the untangling transformation, see Eq.~\eqref{eq:diagvar}, are $\alpha
=-8.46\times 10^{-5}$ and $\beta =2.14\times 10^{-5}\,$.  Thus, the
quadratic part of the Hamiltonian~\eqref{eq:untangling} in diagonal form reads
$$
H_0=2.52\times 10^{-2} \sigma_1^2 +7.00\times 10^{-7}\sigma_3^2
+7.20\times 10^1 \Sigma_1^2 +2.01\times10^2 \Sigma_3^2\ .
$$

The parameters $U_{1}^*$ and $U_{3}^*$ related to the rescaled polar
coordinates, see Eq.~\eqref{eq:actang}, take the values
$U_{1*}=5.348\times 10^{1}$ and $U_{3*}=1.696\times 10^{4}\,$, and the
quadratic part of the Hamiltonian in action angle variables, see
Eq.~\eqref{eq:H_Uu}, reads
$$
H_0=2.69\,U_1+2.37\times 10^{-2}U_3\ ,
$$
while the term of order $3$ in Eq.~\eqref{eq:H_Uu}, $H^{(0)}_1\,$, reads
\begin{align*}
H^{(0)}_1 = &
 2.86\times 10^{-7}\cos (u_1)   \sqrt{U_1}^{3}
-2.87\times 10^{-7}\cos (3 u_1) \sqrt{U_1}^{3} \\
&
+1.98\times 10^{-3}\cos (2u_1-u_3)  U_1 \sqrt{U_3}
-3.99\times 10^{-3}\cos (u_3)       U_1 \sqrt{U_3} \\
&
+2.01\times 10^{-3}\cos (2u_1+u_3) U_1 \sqrt{U_3}
-3.97\times 10^{-3}\cos (u_1)      \sqrt{U_1} U_3 \\
&
+9.29\times 10^{-2}\cos (u_1-2u_3) \sqrt{U_1} U_3
-9.62\times 10^{-2}\cos (u_1+2u_3) \sqrt{U_1} U_3 \\
&
-2.19 \cos (u_3)  \sqrt{U_3}^{3}
-2.19 \cos (3u_3) \sqrt{U_3}^{3}\ .
\end{align*}

We provide the approximation of Eq.~\eqref{eq:H_Uu}, up to order $6$
in $\left(\sqrt{U_1}\,,\ \sqrt{U_3}\,\right)\,$, in the electronic
Supplemental Material (see Table~2) while we report below the number
of coefficients in Eq.~\eqref{eq:H_Uu} at each order,
\begin{center}
\begin{tabular}{ c | c@{\hskip 5pt}c@{\hskip 5pt}c@{\hskip 5pt}c@{\hskip 5pt}c@{\hskip 5pt}c@{\hskip 5pt}c@{\hskip 5pt}c@{\hskip 5pt}c@{\hskip 5pt}c@{\hskip 5pt}c@{\hskip 5pt}c@{\hskip 5pt}c@{\hskip 5pt}c@{\hskip 5pt}c }
Order    & 2&3& 4 &  5 &  6 &  7 &  8 &  9 & 10 &  11 &  12 &  13 &  14 &  15 &  16 \\
\hline
\#terms & 2&10& 19 & 28 & 44 & 54 & 70 & 84 & 93 & 105 & 112 & 125 & 130 & 143 & 145 \\
\end{tabular}\ .
\end{center}

Finally, we report in the electronic Supplemental Material (see
Table~3) the truncated normal form, up to order $12$.  In this case
the number of terms of order $2r$ in
$\left(\sqrt{U_1}\,,\ \sqrt{U_3}\,\right)$ is just equal to $r+1\,$,
thus, at each order, we report in the table below the number of
coefficients in the remainder term at different orders
\begin{center}
\begin{tabular}{ c | c@{\hskip 5pt}c@{\hskip 5pt}c@{\hskip 5pt}c@{\hskip 5pt}c@{\hskip 5pt}c@{\hskip 5pt}c@{\hskip 5pt}c@{\hskip 5pt}c@{\hskip 5pt}c@{\hskip 5pt}c@{\hskip 5pt}c@{\hskip 5pt}c@{\hskip 5pt}c }
Order    & 3  & 4 &  5 &  6 &  7 &  8 &  9 & 10 &  11 &  12 &  13 &  14 &  15 &  16 \\
\hline
\#terms  & 10 & 19 & 28 & 44 & 60 & 85 & 110 & 146 & 182 & 231 & 280 & 345 & 423 & 544 \\
\end{tabular}\ .
\end{center}

\end{document}

%% file: time_rho0_U3.tex
\begingroup
  \makeatletter
  \providecommand\color[2][]{%
    \GenericError{(gnuplot) \space\space\space\@spaces}{%
      Package color not loaded in conjunction with
      terminal option `colourtext'%
    }{See the gnuplot documentation for explanation.%
    }{Either use 'blacktext' in gnuplot or load the package
      color.sty in LaTeX.}%
    \renewcommand\color[2][]{}%
  }%
  \providecommand\includegraphics[2][]{%
    \GenericError{(gnuplot) \space\space\space\@spaces}{%
      Package graphicx or graphics not loaded%
    }{See the gnuplot documentation for explanation.%
    }{The gnuplot epslatex terminal needs graphicx.sty or graphics.sty.}%
    \renewcommand\includegraphics[2][]{}%
  }%
  \providecommand\rotatebox[2]{#2}%
  \@ifundefined{ifGPcolor}{%
    \newif\ifGPcolor
    \GPcolortrue
  }{}%
  \@ifundefined{ifGPblacktext}{%
    \newif\ifGPblacktext
    \GPblacktexttrue
  }{}%
  \let\gplgaddtomacro\g@addto@macro
  \gdef\gplbacktext{}%
  \gdef\gplfronttext{}%
  \makeatother
  \ifGPblacktext
    \def\colorrgb#1{}%
    \def\colorgray#1{}%
  \else
    \ifGPcolor
      \def\colorrgb#1{\color[rgb]{#1}}%
      \def\colorgray#1{\color[gray]{#1}}%
      \expandafter\def\csname LTw\endcsname{\color{white}}%
      \expandafter\def\csname LTb\endcsname{\color{black}}%
      \expandafter\def\csname LTa\endcsname{\color{black}}%
      \expandafter\def\csname LT0\endcsname{\color[rgb]{1,0,0}}%
      \expandafter\def\csname LT1\endcsname{\color[rgb]{0,1,0}}%
      \expandafter\def\csname LT2\endcsname{\color[rgb]{0,0,1}}%
      \expandafter\def\csname LT3\endcsname{\color[rgb]{1,0,1}}%
      \expandafter\def\csname LT4\endcsname{\color[rgb]{0,1,1}}%
      \expandafter\def\csname LT5\endcsname{\color[rgb]{1,1,0}}%
      \expandafter\def\csname LT6\endcsname{\color[rgb]{0,0,0}}%
      \expandafter\def\csname LT7\endcsname{\color[rgb]{1,0.3,0}}%
      \expandafter\def\csname LT8\endcsname{\color[rgb]{0.5,0.5,0.5}}%
    \else
      \def\colorrgb#1{\color{black}}%
      \def\colorgray#1{\color[gray]{#1}}%
      \expandafter\def\csname LTw\endcsname{\color{white}}%
      \expandafter\def\csname LTb\endcsname{\color{black}}%
      \expandafter\def\csname LTa\endcsname{\color{black}}%
      \expandafter\def\csname LT0\endcsname{\color{black}}%
      \expandafter\def\csname LT1\endcsname{\color{black}}%
      \expandafter\def\csname LT2\endcsname{\color{black}}%
      \expandafter\def\csname LT3\endcsname{\color{black}}%
      \expandafter\def\csname LT4\endcsname{\color{black}}%
      \expandafter\def\csname LT5\endcsname{\color{black}}%
      \expandafter\def\csname LT6\endcsname{\color{black}}%
      \expandafter\def\csname LT7\endcsname{\color{black}}%
      \expandafter\def\csname LT8\endcsname{\color{black}}%
    \fi
  \fi
  \setlength{\unitlength}{0.0500bp}%
  \begin{picture}(7200.00,5040.00)%
    \gplgaddtomacro\gplbacktext{%
      \csname LTb\endcsname%
      \put(814,704){\makebox(0,0)[r]{\strut{} 0}}%
      \put(814,1317){\makebox(0,0)[r]{\strut{} 10}}%
      \put(814,1929){\makebox(0,0)[r]{\strut{} 20}}%
      \put(814,2542){\makebox(0,0)[r]{\strut{} 30}}%
      \put(814,3154){\makebox(0,0)[r]{\strut{} 40}}%
      \put(814,3767){\makebox(0,0)[r]{\strut{} 50}}%
      \put(814,4379){\makebox(0,0)[r]{\strut{} 60}}%
      \put(946,484){\makebox(0,0){\strut{} 0}}%
      \put(2117,484){\makebox(0,0){\strut{} 1}}%
      \put(3289,484){\makebox(0,0){\strut{} 2}}%
      \put(4460,484){\makebox(0,0){\strut{} 3}}%
      \put(5632,484){\makebox(0,0){\strut{} 4}}%
      \put(6803,484){\makebox(0,0){\strut{} 5}}%
      \put(176,2541){\rotatebox{-270}{\makebox(0,0){\strut{}$\boldsymbol{\log_{10} T(\rho_0 )}$}}}%
      \put(3874,154){\makebox(0,0){\strut{}$\boldsymbol{\rho_0}$}}%
      \put(3874,4709){\makebox(0,0){\bf Effective stability time}}%
    }%
    \gplgaddtomacro\gplfronttext{%
    }%
    \gplbacktext
    \put(0,0){\includegraphics{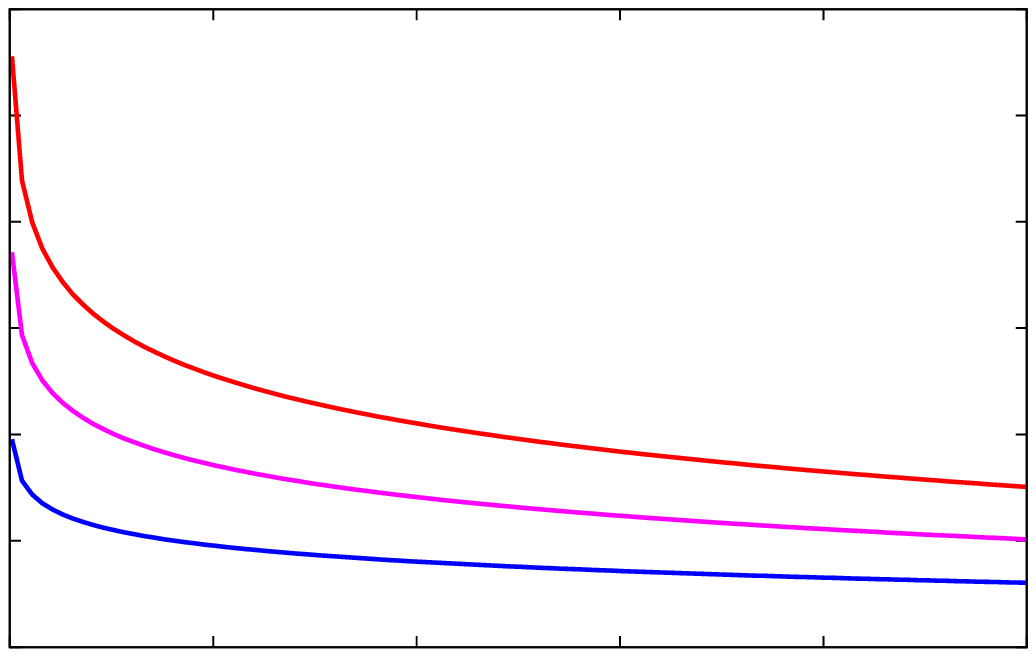}}%
    \gplfronttext
  \end{picture}%
\endgroup

%% file: time_Omdot-incl_U3.tex
\begingroup
  \makeatletter
  \providecommand\color[2][]{%
    \GenericError{(gnuplot) \space\space\space\@spaces}{%
      Package color not loaded in conjunction with
      terminal option `colourtext'%
    }{See the gnuplot documentation for explanation.%
    }{Either use 'blacktext' in gnuplot or load the package
      color.sty in LaTeX.}%
    \renewcommand\color[2][]{}%
  }%
  \providecommand\includegraphics[2][]{%
    \GenericError{(gnuplot) \space\space\space\@spaces}{%
      Package graphicx or graphics not loaded%
    }{See the gnuplot documentation for explanation.%
    }{The gnuplot epslatex terminal needs graphicx.sty or graphics.sty.}%
    \renewcommand\includegraphics[2][]{}%
  }%
  \providecommand\rotatebox[2]{#2}%
  \@ifundefined{ifGPcolor}{%
    \newif\ifGPcolor
    \GPcolortrue
  }{}%
  \@ifundefined{ifGPblacktext}{%
    \newif\ifGPblacktext
    \GPblacktexttrue
  }{}%
  \let\gplgaddtomacro\g@addto@macro
  \gdef\gplbacktext{}%
  \gdef\gplfronttext{}%
  \makeatother
  \ifGPblacktext
    \def\colorrgb#1{}%
    \def\colorgray#1{}%
  \else
    \ifGPcolor
      \def\colorrgb#1{\color[rgb]{#1}}%
      \def\colorgray#1{\color[gray]{#1}}%
      \expandafter\def\csname LTw\endcsname{\color{white}}%
      \expandafter\def\csname LTb\endcsname{\color{black}}%
      \expandafter\def\csname LTa\endcsname{\color{black}}%
      \expandafter\def\csname LT0\endcsname{\color[rgb]{1,0,0}}%
      \expandafter\def\csname LT1\endcsname{\color[rgb]{0,1,0}}%
      \expandafter\def\csname LT2\endcsname{\color[rgb]{0,0,1}}%
      \expandafter\def\csname LT3\endcsname{\color[rgb]{1,0,1}}%
      \expandafter\def\csname LT4\endcsname{\color[rgb]{0,1,1}}%
      \expandafter\def\csname LT5\endcsname{\color[rgb]{1,1,0}}%
      \expandafter\def\csname LT6\endcsname{\color[rgb]{0,0,0}}%
      \expandafter\def\csname LT7\endcsname{\color[rgb]{1,0.3,0}}%
      \expandafter\def\csname LT8\endcsname{\color[rgb]{0.5,0.5,0.5}}%
    \else
      \def\colorrgb#1{\color{black}}%
      \def\colorgray#1{\color[gray]{#1}}%
      \expandafter\def\csname LTw\endcsname{\color{white}}%
      \expandafter\def\csname LTb\endcsname{\color{black}}%
      \expandafter\def\csname LTa\endcsname{\color{black}}%
      \expandafter\def\csname LT0\endcsname{\color{black}}%
      \expandafter\def\csname LT1\endcsname{\color{black}}%
      \expandafter\def\csname LT2\endcsname{\color{black}}%
      \expandafter\def\csname LT3\endcsname{\color{black}}%
      \expandafter\def\csname LT4\endcsname{\color{black}}%
      \expandafter\def\csname LT5\endcsname{\color{black}}%
      \expandafter\def\csname LT6\endcsname{\color{black}}%
      \expandafter\def\csname LT7\endcsname{\color{black}}%
      \expandafter\def\csname LT8\endcsname{\color{black}}%
    \fi
  \fi
  \setlength{\unitlength}{0.0500bp}%
  \begin{picture}(7200.00,5040.00)%
    \gplgaddtomacro\gplbacktext{%
      \csname LTb\endcsname%
      \put(3600,4312){\makebox(0,0){\bf Effective stability time ($\boldsymbol{\dot\Omega}$ vs. $\boldsymbol{i}$)}}%
    }%
    \gplgaddtomacro\gplfronttext{%
      \csname LTb\endcsname%
      \put(5355,772){\makebox(0,0){\strut{}-0.016}}%
      \put(4629,772){\makebox(0,0){\strut{}-0.014}}%
      \put(3903,772){\makebox(0,0){\strut{}-0.012}}%
      \put(3178,772){\makebox(0,0){\strut{}-0.010}}%
      \put(2452,772){\makebox(0,0){\strut{}-0.008}}%
      \put(1726,772){\makebox(0,0){\strut{}-0.006}}%
      \put(3600,442){\makebox(0,0){\strut{}$\boldsymbol{\dot\Omega}$}}%
      \put(998,1414){\makebox(0,0)[r]{\strut{}0.004}}%
      \put(998,2221){\makebox(0,0)[r]{\strut{}0.008}}%
      \put(998,3026){\makebox(0,0)[r]{\strut{}0.012}}%
      \put(998,3833){\makebox(0,0)[r]{\strut{}0.016}}%
      \put(272,2520){\rotatebox{-270}{\makebox(0,0){\strut{}$\boldsymbol{i}$}}}%
      \put(6527,1058){\makebox(0,0)[l]{\strut{} 12}}%
      \put(6527,1382){\makebox(0,0)[l]{\strut{} 14}}%
      \put(6527,1707){\makebox(0,0)[l]{\strut{} 16}}%
      \put(6527,2032){\makebox(0,0)[l]{\strut{} 18}}%
      \put(6527,2357){\makebox(0,0)[l]{\strut{} 20}}%
      \put(6527,2682){\makebox(0,0)[l]{\strut{} 22}}%
      \put(6527,3007){\makebox(0,0)[l]{\strut{} 24}}%
      \put(6527,3332){\makebox(0,0)[l]{\strut{} 26}}%
      \put(6527,3657){\makebox(0,0)[l]{\strut{} 28}}%
      \put(6527,3982){\makebox(0,0)[l]{\strut{} 30}}%
      \put(6950,2520){\rotatebox{-270}{\makebox(0,0){\strut{}$\boldsymbol{\log_{10}T}$}}}%
    }%
    \gplbacktext
    \put(0,0){\includegraphics{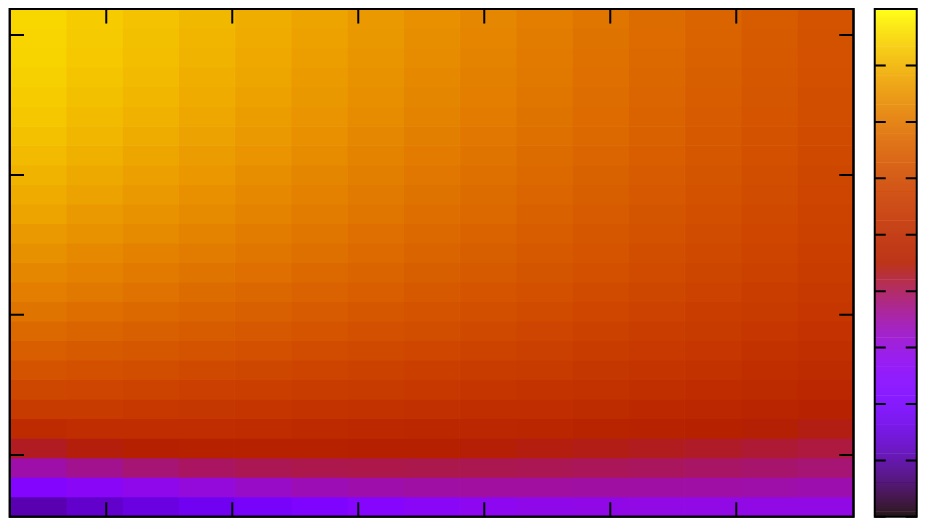}}%
    \gplfronttext
  \end{picture}%
\endgroup

%% file: time_Omega-CRM_U3.tex
\begingroup
  \makeatletter
  \providecommand\color[2][]{%
    \GenericError{(gnuplot) \space\space\space\@spaces}{%
      Package color not loaded in conjunction with
      terminal option `colourtext'%
    }{See the gnuplot documentation for explanation.%
    }{Either use 'blacktext' in gnuplot or load the package
      color.sty in LaTeX.}%
    \renewcommand\color[2][]{}%
  }%
  \providecommand\includegraphics[2][]{%
    \GenericError{(gnuplot) \space\space\space\@spaces}{%
      Package graphicx or graphics not loaded%
    }{See the gnuplot documentation for explanation.%
    }{The gnuplot epslatex terminal needs graphicx.sty or graphics.sty.}%
    \renewcommand\includegraphics[2][]{}%
  }%
  \providecommand\rotatebox[2]{#2}%
  \@ifundefined{ifGPcolor}{%
    \newif\ifGPcolor
    \GPcolortrue
  }{}%
  \@ifundefined{ifGPblacktext}{%
    \newif\ifGPblacktext
    \GPblacktexttrue
  }{}%
  \let\gplgaddtomacro\g@addto@macro
  \gdef\gplbacktext{}%
  \gdef\gplfronttext{}%
  \makeatother
  \ifGPblacktext
    \def\colorrgb#1{}%
    \def\colorgray#1{}%
  \else
    \ifGPcolor
      \def\colorrgb#1{\color[rgb]{#1}}%
      \def\colorgray#1{\color[gray]{#1}}%
      \expandafter\def\csname LTw\endcsname{\color{white}}%
      \expandafter\def\csname LTb\endcsname{\color{black}}%
      \expandafter\def\csname LTa\endcsname{\color{black}}%
      \expandafter\def\csname LT0\endcsname{\color[rgb]{1,0,0}}%
      \expandafter\def\csname LT1\endcsname{\color[rgb]{0,1,0}}%
      \expandafter\def\csname LT2\endcsname{\color[rgb]{0,0,1}}%
      \expandafter\def\csname LT3\endcsname{\color[rgb]{1,0,1}}%
      \expandafter\def\csname LT4\endcsname{\color[rgb]{0,1,1}}%
      \expandafter\def\csname LT5\endcsname{\color[rgb]{1,1,0}}%
      \expandafter\def\csname LT6\endcsname{\color[rgb]{0,0,0}}%
      \expandafter\def\csname LT7\endcsname{\color[rgb]{1,0.3,0}}%
      \expandafter\def\csname LT8\endcsname{\color[rgb]{0.5,0.5,0.5}}%
    \else
      \def\colorrgb#1{\color{black}}%
      \def\colorgray#1{\color[gray]{#1}}%
      \expandafter\def\csname LTw\endcsname{\color{white}}%
      \expandafter\def\csname LTb\endcsname{\color{black}}%
      \expandafter\def\csname LTa\endcsname{\color{black}}%
      \expandafter\def\csname LT0\endcsname{\color{black}}%
      \expandafter\def\csname LT1\endcsname{\color{black}}%
      \expandafter\def\csname LT2\endcsname{\color{black}}%
      \expandafter\def\csname LT3\endcsname{\color{black}}%
      \expandafter\def\csname LT4\endcsname{\color{black}}%
      \expandafter\def\csname LT5\endcsname{\color{black}}%
      \expandafter\def\csname LT6\endcsname{\color{black}}%
      \expandafter\def\csname LT7\endcsname{\color{black}}%
      \expandafter\def\csname LT8\endcsname{\color{black}}%
    \fi
  \fi
  \setlength{\unitlength}{0.0500bp}%
  \begin{picture}(7200.00,5040.00)%
    \gplgaddtomacro\gplbacktext{%
      \csname LTb\endcsname%
      \put(3600,4312){\makebox(0,0){\bf Effective stability time ($\boldsymbol{\dot\Omega}$ vs. $\boldsymbol{C/MR_{\rm e}^2}$)}}%
    }%
    \gplgaddtomacro\gplfronttext{%
      \csname LTb\endcsname%
      \put(5355,772){\makebox(0,0){\strut{}-0.016}}%
      \put(4629,772){\makebox(0,0){\strut{}-0.014}}%
      \put(3903,772){\makebox(0,0){\strut{}-0.012}}%
      \put(3178,772){\makebox(0,0){\strut{}-0.010}}%
      \put(2452,772){\makebox(0,0){\strut{}-0.008}}%
      \put(1726,772){\makebox(0,0){\strut{}-0.006}}%
      \put(3600,442){\makebox(0,0){\strut{}$\boldsymbol{\dot\Omega}$}}%
      \put(998,1058){\makebox(0,0)[r]{\strut{}0.30}}%
      \put(998,1643){\makebox(0,0)[r]{\strut{}0.32}}%
      \put(998,2228){\makebox(0,0)[r]{\strut{}0.34}}%
      \put(998,2812){\makebox(0,0)[r]{\strut{}0.36}}%
      \put(998,3397){\makebox(0,0)[r]{\strut{}0.38}}%
      \put(998,3982){\makebox(0,0)[r]{\strut{}0.40}}%
      \put(404,2520){\rotatebox{-270}{\makebox(0,0){\strut{}$\boldsymbol{C/MR_{\rm e}^2}$}}}%
      \put(6527,1058){\makebox(0,0)[l]{\strut{} 20}}%
      \put(6527,1423){\makebox(0,0)[l]{\strut{} 20.5}}%
      \put(6527,1789){\makebox(0,0)[l]{\strut{} 21}}%
      \put(6527,2154){\makebox(0,0)[l]{\strut{} 21.5}}%
      \put(6527,2520){\makebox(0,0)[l]{\strut{} 22}}%
      \put(6527,2885){\makebox(0,0)[l]{\strut{} 22.5}}%
      \put(6527,3251){\makebox(0,0)[l]{\strut{} 23}}%
      \put(6527,3616){\makebox(0,0)[l]{\strut{} 23.5}}%
      \put(6527,3982){\makebox(0,0)[l]{\strut{} 24}}%
      \put(7050,2520){\rotatebox{-270}{\makebox(0,0){\strut{}$\boldsymbol{\log_{10}T}$}}}%
    }%
    \gplbacktext
    \put(0,0){\includegraphics{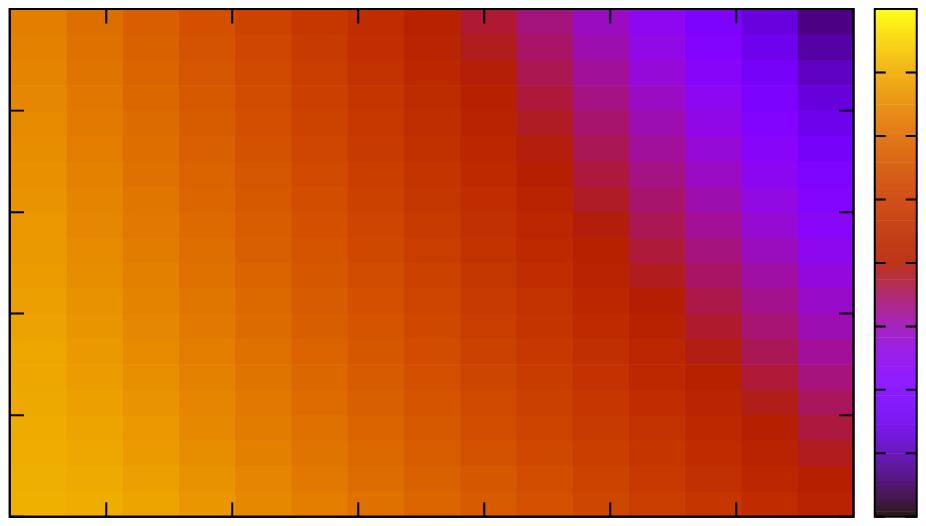}}%
    \gplfronttext
  \end{picture}%
\endgroup

%% file: time_CRM-incl_U3.tex
\begingroup
  \makeatletter
  \providecommand\color[2][]{%
    \GenericError{(gnuplot) \space\space\space\@spaces}{%
      Package color not loaded in conjunction with
      terminal option `colourtext'%
    }{See the gnuplot documentation for explanation.%
    }{Either use 'blacktext' in gnuplot or load the package
      color.sty in LaTeX.}%
    \renewcommand\color[2][]{}%
  }%
  \providecommand\includegraphics[2][]{%
    \GenericError{(gnuplot) \space\space\space\@spaces}{%
      Package graphicx or graphics not loaded%
    }{See the gnuplot documentation for explanation.%
    }{The gnuplot epslatex terminal needs graphicx.sty or graphics.sty.}%
    \renewcommand\includegraphics[2][]{}%
  }%
  \providecommand\rotatebox[2]{#2}%
  \@ifundefined{ifGPcolor}{%
    \newif\ifGPcolor
    \GPcolortrue
  }{}%
  \@ifundefined{ifGPblacktext}{%
    \newif\ifGPblacktext
    \GPblacktexttrue
  }{}%
  \let\gplgaddtomacro\g@addto@macro
  \gdef\gplbacktext{}%
  \gdef\gplfronttext{}%
  \makeatother
  \ifGPblacktext
    \def\colorrgb#1{}%
    \def\colorgray#1{}%
  \else
    \ifGPcolor
      \def\colorrgb#1{\color[rgb]{#1}}%
      \def\colorgray#1{\color[gray]{#1}}%
      \expandafter\def\csname LTw\endcsname{\color{white}}%
      \expandafter\def\csname LTb\endcsname{\color{black}}%
      \expandafter\def\csname LTa\endcsname{\color{black}}%
      \expandafter\def\csname LT0\endcsname{\color[rgb]{1,0,0}}%
      \expandafter\def\csname LT1\endcsname{\color[rgb]{0,1,0}}%
      \expandafter\def\csname LT2\endcsname{\color[rgb]{0,0,1}}%
      \expandafter\def\csname LT3\endcsname{\color[rgb]{1,0,1}}%
      \expandafter\def\csname LT4\endcsname{\color[rgb]{0,1,1}}%
      \expandafter\def\csname LT5\endcsname{\color[rgb]{1,1,0}}%
      \expandafter\def\csname LT6\endcsname{\color[rgb]{0,0,0}}%
      \expandafter\def\csname LT7\endcsname{\color[rgb]{1,0.3,0}}%
      \expandafter\def\csname LT8\endcsname{\color[rgb]{0.5,0.5,0.5}}%
    \else
      \def\colorrgb#1{\color{black}}%
      \def\colorgray#1{\color[gray]{#1}}%
      \expandafter\def\csname LTw\endcsname{\color{white}}%
      \expandafter\def\csname LTb\endcsname{\color{black}}%
      \expandafter\def\csname LTa\endcsname{\color{black}}%
      \expandafter\def\csname LT0\endcsname{\color{black}}%
      \expandafter\def\csname LT1\endcsname{\color{black}}%
      \expandafter\def\csname LT2\endcsname{\color{black}}%
      \expandafter\def\csname LT3\endcsname{\color{black}}%
      \expandafter\def\csname LT4\endcsname{\color{black}}%
      \expandafter\def\csname LT5\endcsname{\color{black}}%
      \expandafter\def\csname LT6\endcsname{\color{black}}%
      \expandafter\def\csname LT7\endcsname{\color{black}}%
      \expandafter\def\csname LT8\endcsname{\color{black}}%
    \fi
  \fi
  \setlength{\unitlength}{0.0500bp}%
  \begin{picture}(7200.00,5040.00)%
    \gplgaddtomacro\gplbacktext{%
      \csname LTb\endcsname%
      \put(3600,4312){\makebox(0,0){\bf Effective stability time ($\boldsymbol{i}$ vs. $\boldsymbol{C/MR_{\rm e}^2}$)}}%
    }%
    \gplgaddtomacro\gplfronttext{%
      \csname LTb\endcsname%
      \put(1762,772){\makebox(0,0){\strut{}0.004}}%
      \put(3103,772){\makebox(0,0){\strut{}0.008}}%
      \put(4442,772){\makebox(0,0){\strut{}0.012}}%
      \put(5783,772){\makebox(0,0){\strut{}0.016}}%
      \put(3600,442){\makebox(0,0){\strut{}$\boldsymbol{i}$}}%
      \put(998,1058){\makebox(0,0)[r]{\strut{}0.30}}%
      \put(998,1643){\makebox(0,0)[r]{\strut{}0.32}}%
      \put(998,2228){\makebox(0,0)[r]{\strut{}0.34}}%
      \put(998,2812){\makebox(0,0)[r]{\strut{}0.36}}%
      \put(998,3397){\makebox(0,0)[r]{\strut{}0.38}}%
      \put(998,3982){\makebox(0,0)[r]{\strut{}0.40}}%
      \put(404,2520){\rotatebox{-270}{\makebox(0,0){\strut{}$\boldsymbol{C/MR_{\rm e}^2}$}}}%
      \put(6527,1058){\makebox(0,0)[l]{\strut{} 14}}%
      \put(6527,1423){\makebox(0,0)[l]{\strut{} 16}}%
      \put(6527,1789){\makebox(0,0)[l]{\strut{} 18}}%
      \put(6527,2154){\makebox(0,0)[l]{\strut{} 20}}%
      \put(6527,2520){\makebox(0,0)[l]{\strut{} 22}}%
      \put(6527,2885){\makebox(0,0)[l]{\strut{} 24}}%
      \put(6527,3251){\makebox(0,0)[l]{\strut{} 26}}%
      \put(6527,3616){\makebox(0,0)[l]{\strut{} 28}}%
      \put(6527,3982){\makebox(0,0)[l]{\strut{} 30}}%
      \put(6950,2520){\rotatebox{-270}{\makebox(0,0){\strut{}$\boldsymbol{\log_{10} T}$}}}%
    }%
    \gplbacktext
    \put(0,0){\includegraphics{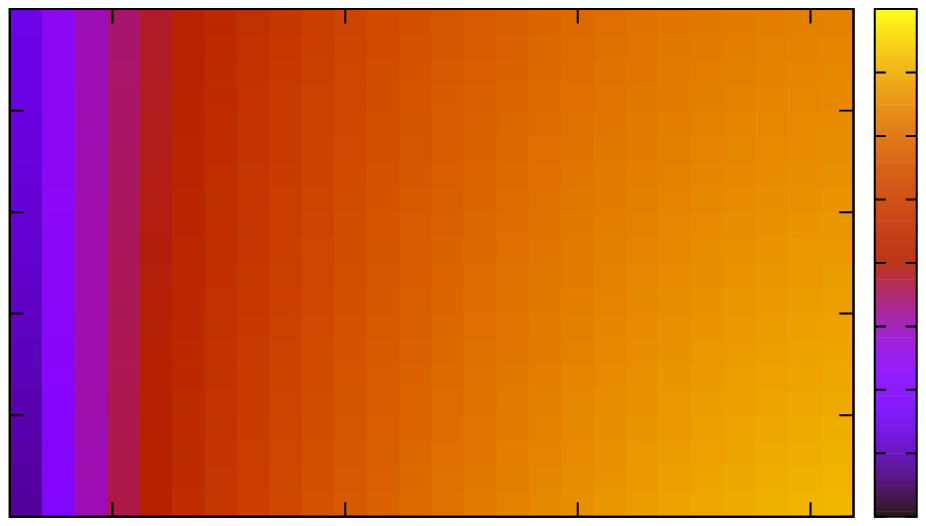}}%
    \gplfronttext
  \end{picture}%
\endgroup